\documentclass[12pt]{article}

\usepackage[T1]{fontenc}
\usepackage{avant}
\usepackage[cp1250]{inputenc}
\usepackage{amsmath,amssymb}
\usepackage{graphicx}
\usepackage{color}
\usepackage{amsmath}
\usepackage{xcolor}

\newtheorem{theorem}{Theorem}
\newtheorem{proposition}{Proposition}
\newtheorem{lemma}{Lemma}
\newtheorem{definition}{Definition}

\newtheorem{corollary}{Corollary}

\newtheorem{question}{Question}
\newtheorem{conjecture}{Conjecture}

\title{On locally homogeneous pseudo-Riemannian compact Einstein manifolds}
\author{Maciej Boche\'nski and Aleksy Tralle}

\begin{document}

\maketitle{}
\abstract{We ask a general question: what are locally homogeneous compact pseudo-Riemannian Einstein manifolds? We show that any standard compact Clifford-Klein form of a simple non-compact Lie group admits at least one Einstein metric. We conjecture that these are basically the only possible locally homogeneous Einstein Riemannian compact manifolds using T. Kobayashi's conjecture as a guiding principle. }
\vskip6pt
\section{Introduction}
Let $M$ be a compact pseudo-Riemannian manifold. We say that $M$ is locally homogeneous if   $M$ is diffeomorphic to the double quotient $\Gamma\backslash G/H$ determined by a discrete subgroup $\Gamma\subset G$ acting freely and properly on $G/H$.  Thus, $M$ admits a system of local coordinates modeled on a homogeneous space $G/H$. For example, one can consider physically interesting manifolds locally modeled on anti-de Sitter space $AdS^n=PO(n-1,2)/O(n-1,1)$, or locally hyperbolic manifolds modeled on the hyperbolic homogeneous space $\mathbb{H}^n=PO(n,1)/O(n)$. Note that  $\Gamma\backslash G/H$ may be  compact in some cases, but it may also happen that no $\Gamma$ can yield a compact locally homogeneous manifold. For example, Calabi and Markus \cite{CM} considered a relativistic spherical space form problem and classified all such forms satisfying the perfect cosmological principle. In particular, they showed that a relativistic spherical space form is always non-compact and has a finite fundamental group. Also, they noted that in  such worlds the Einstein field equation is satisfied under some conditions.   
Motivated by this as well as by purely mathematical works on compact space forms \cite{kas}, \cite{kob}, \cite{koy},\cite{ku} we ask the following more general question.
\begin{question} What are locally homogeneous compact pseudo-Riemannian Einstein manifolds?
\end{question}

It is well known that mathematicians tend to consider the Einstein equation 
$\operatorname{Ric}(\rho)=\lambda\cdot \rho,\,\lambda=\operatorname{const}$
on compact  {\it Riemannian manifolds}, while physicists want to solve this equation in the compact {\it pseudo-Riemannian} case. However, there are various obstructions to even put a pseudo-Riemannian metric on a compact manifold, and adding the Einstein equation might yield  even more difficulties. This fact is one of the main motivations of this work.  Let $G/H$ be a {\it non-compact} homogeneous space of a semisimple connected Lie group $G$ and a {\it non-compact} closed subgroup $H$. This homogeneous space carries a $G$-invariant pseudo-Riemannian metrics (for example, there is at least one induced by the Killing form). Assume that there is a discrete subgroup $\Gamma\subset G$ which acts freely and properly on $G/H$ with a compact quotient. Then one obtains the compact manifold $\Gamma\setminus G/H$ which inherits the pseudo-Riemannian metric.  Note that if $H$ is non-compact, the existence of such $\Gamma$ becomes a problem very difficult to solve. In the case of compact $H$ the solution is much easier: one just takes a co-compact lattice $\Gamma$ in $G$. It is easy to see that if $H$ is compact, $\Gamma\setminus G/H$ is a compact orbifold, and by Selberg's lemma, one can find a finite index subgroup $\Gamma'\subset \Gamma$ which acts freely on $G/H$. All this is not possible for non-compact $H$. However, T. Kobayashi \cite{kob} have found a kind of analogue of this. Assume that we are given a homogeneous space $G/H$ of a semisimple connected linear Lie group $G$ and a proper and {\it transitive} action of a reductive subgroup $L\subset G$. Then one obtains a diffeomorphism 
$$L/L\cap H\cong G/H.$$
Assume that $L\cap H$ is compact and let $\Gamma\subset L$ be a co-compact lattice in $L$. One obtains
$$\Gamma\backslash G/H\cong \Gamma\backslash L/L\cap H.$$
Since $L\cap H$ is compact and $\Gamma$ is a co-compact lattice in $L$, $\Gamma\backslash L/L\cap H$ is a compact orbifold, and the same is valid for $\Gamma\backslash G/H$. Again one can use Selberg's Lemma and take $\Gamma ' \subset \Gamma$ so that $\Gamma ' \backslash G/H$ is a compact manifold.
\begin{definition}
 {\rm We say that $\Gamma\backslash G/H$ is a (compact) Clifford-Klein form, if  $\Gamma$ is a discrete subgroup of $G$ acting properly (and co-compactly) on $G/H$.}
 \end{definition}  
 
\begin{definition} {\rm Let $G$ be a semisimple non-compact, linear Lie group and let   $H$ and $L$ be reductive subgroups in $G$. We say that the triple $(G,H,L)$ is an {\it Onishchik triple}, if $G=HL$.
 Let $(G,H,L)$ be an Onishchik triple. If $H\cap L$ is compact, we say that $(G,H,L)$ is a {\it standard triple}. We say that $\Gamma\backslash G/H$ is a {\it standard Clifford-Klein form} if it is obtained from some standard triple $(G,H,L)$.}

\end{definition}
The name ``Onishchik triples'' is justified by the fact that they were first considered and classified by A. L. Onishchik \cite{o1}.
\begin{conjecture}[T. Kobayashi, \cite{koy}] If a homogeneous space $G/H$ of a connected semisimple non-compact Lie group $G$ has a compact Clifford-Klein form, then there is a structure of a Clifford-Klein form on $G/H$ which is standard.
\end{conjecture}
 Note that this conjecture does not say that all compact Clifford-Klein forms are standard, it says that for $\Gamma\backslash G/H$ there should exist a standard triple $(G,H,\hat L)$. {\it All known examples of Clifford-Klein forms are either standard, or are determined by some standard triples $(G,H,L)$ but with a discrete subgroup $\hat\Gamma$ which is a deformation of a subgroup $\Gamma\subset L$} (see \cite{kas},\cite{koy}). The deformation can be chosen in a way that $\hat\Gamma$ is Zariski dense in $G$ \cite{kas}. 
\vskip6pt
There are many partial results which yield a strong evidence for the conjecture, for example \cite{bt},\cite{koy},\cite{mor},\cite{th}.  However, the final answer (which is challenging and mathematically interesting) is still to be answered. On the other hand, in view of this we are looking for locally homogeneous Einstein manifolds using the Kobayashi conjecture as a guiding principle: {\it what are possible Einstein Clifford-Klein forms?}

The main result of this work is the following.
\begin{theorem}\label{thm:main} Any standard triple $(G,H,L)$ of a simple non-compact Lie group determines a compact Clifford-Klein form with at least one Einstein metric.
\end{theorem}

This leads us to the following hypothesis.
\begin{conjecture}\label{conj:einstein} The only compact locally homogeneous simple pseudo Riemannian Einstein manifolds are of the form
 $\Gamma\backslash G/H$,determined by
 \begin{itemize}
 \item  some standard triple $(G,H,L)$ and a co-compact lattice $\Gamma$ in $L$, or the deformation of such $\Gamma$; 

\item $(G,\hat H,\hat L)$, where $\hat H,\hat L$ differ from $H$ and $L$ by a compact factor.
\end{itemize}
\end{conjecture}
Note that if a homogeneous space $G/H$ has a compact Clifford-Klein form $\Gamma\backslash G/H$ then the same is valid for any homogeneous space $G/\hat H$, where $\hat H$ differs from $H$ by a compact factor. It follows that
 the positive answer to the Kobayashi conjecture will  be a very strong evidence for Conjecture \ref{conj:einstein}. 
 
  It should be mentioned that in this paper we don't propose  new methods of looking for  $G$-invariant Einstein metrics.  Instead, we follow our guiding principle of looking for Einstein metrics on Clifford-Klein forms. Our strategy is to reduce our problem to finding such metrics on  {\it compact} homogeneous spaces, which are dual to non-compact homogeneous spaces $G/H$,    and then applying  the classical methods such as  \cite{wz}. The structure of the article is as follows. We use the duality between homogeneous spaces of non-compact semisimple Lie groups and homogeneous spaces of compact Lie groups and the relation between Einstein metrics on dual spaces.  Then we use the classification of Onishchik triples and eliminate the known cases (Section 3). The remaining cases are analyzed using the approach via Einstein metrics with additional symmetries (Proposition 4, Section 4). We observe that  the cases we need to analyze are covered by Proposition 4 (see Theorem 3, Section 5). 
  
  \noindent {\bf Acknowledgment.} The first named author was supported by the National Science Center (Poland), grant no. 2018/31/D/ST1/00083. The second named author was supported by the National Science Center (Poland), grant no. 2018/31/B/ST1/00053. Our special thanks go to the anonimous referee for valuable advice.

\section{Preliminaries}
In this paper we use the basics of Lie theory without further explanations referring to \cite{ov}. We denote the Lie algebras of the Lie groups $G,H,...$ by the corresponding Gothic letters $\mathfrak{g},\mathfrak{h},...$. 
Recall the notion of the {\it Casimir operator} \cite{be}. Let $q$ be a non-degenerate symmetric $\operatorname{ad}(\mathfrak{k})$-invariant bilinear form on a Lie algebra $\mathfrak{k}$, and $\mu: \mathfrak{k}\rightarrow \mathfrak{gl}(N)$ be an $N$-dimensional linear representation. A Casimir operator is an element $C_{\mu,q}\in\mathfrak{gl}(N)$ given by the formula
$$C_{\mu,q}=-\sum_{j=1}^k\mu(U_j)\circ\mu(V_j),$$
where $U_1,..., U_k$ and $V_1,...,V_k$ denote the dual bases of $\mathfrak{k}$ with respect to $q$. In particular, assume that we are given a semisimple Lie algebra $\mathfrak{g}$ with the Killing form $B$  and a subalgebra $\mathfrak{k}$ with the property that $B|_{\mathfrak{k}}$ is non-degenerate. Consider the $\operatorname{ad}(\mathfrak{k})$-invariant decomposition determined by $B$ 
$$\mathfrak{g}=\mathfrak{k}+\mathfrak{m}.$$
We obtain an isotropy representation $\chi:\mathfrak{k}\rightarrow\mathfrak{gl}(\mathfrak{m})$ and the Casimir operator $C_{\chi, -B}$. One can notice that the restriction of $C_{\chi,-B}$ on $\chi(\mathfrak{k})$-invariant subspaces $\mathfrak{m}_i$ of $\mathfrak{m}$ have the form
$$C_{\chi,-B}|_{\mathfrak{m}_i}=c_i^*\operatorname{Id}_{\mathfrak{m}_i}.$$
The scalars $c_i^*$ are called the {\it Casimir constants}.

A Riemannian metric $\rho$ on a manifold $M$ is called Einstein if $\operatorname{Ric}(\rho)=\lambda\cdot\rho$ for some constant $\lambda\in\mathbb{R}$. Our basic reference on Einstein metrics is \cite{be}. The Einstein metrics of volume $1$ on a compact manifold $M$ are the critical points of the total scalar curvature functional 
$$S(\mathfrak{\rho})=\int_MS(\rho)d\operatorname{vol}_{\rho}$$
on the space $\mathcal{M}_1$ of Riemannian metrics of volume $1$. Note that in this and only in this section we denote by $G/H$ homogeneous spaces of {\it compact} Lie groups. It is known \cite{wz} that on  homogeneous spaces of a compact Lie group $G$ the $G$-invariant Einstein metrics are precisely the critical points of the restriction of $S$ onto the space $\mathcal{M}^G$ of $G$-invariant Riemannian metrics. Therefore, one can look for $G$-invariant Einstein metrics using the  following scheme.  Assume that $\mathfrak{g}$ is semisimple and of compact type, and let $B$ be the {\it negative} of the Killing form on $\mathfrak{g}$. Consider the $B$-orthogonal decomposition 
$$\mathfrak{g}=\mathfrak{h}+\mathfrak{m}.$$
Then the set of $G$-invariant metrics on $G/H$ can be identified with the set of $\operatorname{Ad}(H)$-invariant inner products $\langle\cdot,\cdot\rangle$ on $\mathfrak{m}$. In this case one can express the scalar curvature of any $G$-invariant metric by the formula (see \cite{be}, 7.39):
$$S=\frac{1}{2}\sum_{\alpha}B(e_{\alpha},e_{\alpha})-\frac{1}{4}\sum_{\alpha,\beta}\langle [e_{\alpha},e_{\beta}]_{\mathfrak{m}},[e_{\alpha},e_{\beta}]_{\mathfrak{m}}\rangle, \eqno (1)$$
where $e_{\alpha}$ constitute an orthonormal basis of $\mathfrak{m}$ with respect to $\langle\cdot,\cdot\rangle$, and $[\cdot,\cdot]_{\mathfrak{m}}$ denotes the $\mathfrak{m}$-component of the corresponding Lie bracket. Assume now that we are given {\it any} $B$-orthogonal decomposition
$$\mathfrak{m}=\mathfrak{m}_1+...+\mathfrak{m}_s$$
into $\operatorname{Ad}(H)$-invariant components.  If the $\operatorname{Ad}(H)$-invariant inner product is given by the formula
$$\langle\cdot,\cdot\rangle=x_1B|_{\mathfrak{m}_1}+\cdots x_sB|_{\mathfrak{m}_s},x_i>0\eqno (2)$$
then applying formula (1) to (2) one obtains the expression 
$$S=\frac{1}{2}\sum_{i=1}^s\frac{d_i}{x_i}-\frac{1}{4}\sum_{i,j,k}[k,i,j]\frac{x_k}{x_ix_j}\eqno (3)$$
where $[k,j,j]=\sum(A^{\gamma}_{\alpha\beta})^2$ and $A^{\gamma}_{\alpha\beta}$ are determined by
$$A^{\gamma}_{\beta\alpha}=B([e_{\alpha},e_{\beta}],e_{\gamma}),$$ 
the sum is taken over all indices $\alpha,\beta,\gamma$ with $e_{\alpha}\in\mathfrak{m}_i,e_{\beta}\in\mathfrak{m}_j$ and $e_{\gamma}\in\mathfrak{m}_k$. Here $d_i=\dim\,\mathfrak{m}_i$. Note that $e_{\alpha}$ constitute a $B$-orthogonal basis of $\mathfrak{m}$ adapted to the decomposition, that is, $e_{\alpha}\in\mathfrak{m}_i$ for some $i$ and $\alpha<\beta$ if $i<j$. 
It is important to note that formula $(3)$ is valid {\it for any} decomposition $(2)$. In \cite{wz} it is required that $\mathfrak{m}_i$ should be {\it irreducible} $\operatorname{Ad}(H)$-modules. This is because one needs to consider {\it all} $\operatorname{Ad}(H)$-invariant inner products on $\mathfrak{m}$ in order to find critical points of the functional $S|_{\mathcal{M}^G}$. Clearly, formula $(3)$ yields all $Ad(H)$-invariant inner products on $\mathfrak{m}$, if all $\mathfrak{m}_i$ are $\operatorname{Ad}(H)$-irreducible and inequivalent as representations.  
In this paper we will need the following particular case given by Lemma \ref{lagrange}.  
\textcolor{blue}{This result appeared in several previous works, one of the first was \cite{dyk}.}
\begin{lemma}
 Assume that we are given inclusions of closed subgroups $H\subset K\subset G$ determining a fibration $G/H\rightarrow G/K$ with fiber $K/H$.
Denote by 
$$\mathfrak{g}=\mathfrak{h}+\mathfrak{m}_{1}+\mathfrak{m}_{2}=\mathfrak{k}+\mathfrak{m}_{2}$$
the $B$-orthogonal, $\operatorname{Ad}(H)$-invariant decomposition of $\mathfrak{g}$ where $\mathfrak{k}=\mathfrak{h}+\mathfrak{m}_{1}.$ Let $x_1,x_2 >0$ and set
$$\langle\cdot,\cdot\rangle:=x_{1}B|_{\mathfrak{m}_{1}}+x_{2}B|_{\mathfrak{m}_{2}}. \eqno (4)$$
Let $S$ be the restriction of the scalar curvature functional restricted onto the space of all $G$-invariant Riemannian metrics given by (4) of \textcolor{blue}{volume 1.}  
 The critical points of the functional $S$ correspond to the real roots of the  polynomial
$$[1,2,2](2d_{1}+d_{2})t^{2}+$$
$$+(d_{1}[2,2,2]-2d_{1}d_{2})t+2d_{1}d_{2}-2d_{2}[1,2,2]-d_{2}[1,1,1]. \eqno (5)$$
\label{lagrange}
\end{lemma}

\noindent {\it Proof}.
 Note that $A_{\alpha \beta}^{\gamma}$ is symmetric in all three entries. We see that $\{ \frac{e_{\alpha}}{\sqrt{x_{i}}} \ | \ e_{\alpha}\in \mathfrak{m}_{i} \}$ is an orthonormal basis for $\mathfrak{m}$ with respect to  $\langle\cdot,\cdot\rangle$. Then the equation (3) takes the form 
$$S=\frac{d_{1}}{2x_{1}}+\frac{d_{2}}{2x_{2}}-\frac{[1,1,1]}{4x_{1}}-\frac{[2,2,2]}{4x_{2}}-\frac{[1,2,2]}{4}(\frac{x_{1}}{x_{2}^{2}}+\frac{2}{x_{1}}),$$
where $d_{i}=\dim\,\mathfrak{m}_{i}.$ 
Note that we have used the condition $[\mathfrak{m}_{1},\mathfrak{m}_{1}]\subset \mathfrak{k}$, which implies  $[1,1,2]=0$.
Using the Lagrange multipliers method with the volume condition $x_{1}^{d_{1}}x_{2}^{d_{2}}=c=\dim\,\mathfrak{g}$ and setting $t:=\frac{x_{1}}{x_{2}}$ one completes the proof.

\hfill$\square$

In our proofs we will also need a relationship between $[i,j,k]$ and the Casimir operators $C_{\chi,B}|_{\mathfrak{m}_i}$ established in \cite{wz}.
\begin{lemma}[\cite{wz}, Lemma 1.5]\label{lemma:casimir} The following formula relates the numbers $[i,j,k]$ and the Casimir constants:
$$\sum_{j,k}[i,j,k]=d_i(1-2c_i^*).$$
\end{lemma}

\section{Ingredients of proof of Theorem \ref{thm:main}}
In the proof of Theorem \ref{thm:main} we use the duality between homogeneous spaces of semisimple Lie groups of compact and non-compact type,  and the classification of Onishchik triples. The method of the proof of the main theorem of this article is as follows. We will use a theorem of I. Kath (see Theorem \ref{thm:kath}) which shows that  to classify all invariant Einstein metrics on $G/H$ of non-compact Lie group $G$, it is sufficient to classify all invariant Einstein metrics on the compact dual $G_{u}/H_{u}$. Proposition \ref{prop:ko-spaces} shows that we need to analyze compact homogeneous spaces from Table 2 (which we call Kobayashi-Onishchik triples).  
\subsection{Duality}
We consider homogeneous spaces $G/H$  and assume that $G$ is connected linear semisimple, and that $H$ is reductive in $G$.   Under these assumptions every $G/H$ has a dual $G_u/H_u$, which we now describe (see \cite{kobo}, Section 3).  
   Let $G_{u}$ be a compact real form of a (connected) complexification $G^{\mathbb{C}}$ of $G$ and let $H_{u}$ be a compact real form of $H^{\mathbb{C}}\subset G^{\mathbb{C}}$. The space $G_{u}/H_{u}$ is called the  homogeneous space of compact type associated with $G/H$ (or dual to $G/H$). Groups $G_{u},$ $H_{u}$ are called the compact duals of $G$ and   $H,$ respectively.
\subsection{Onishchik triples and standard triples}
Recall that $(G,H,L)$ is an Onishchik triple, if $G=HL$ (that is, any element in $G$ can be represented as a product of an element in $H$ and an element in $L$). It is easy to show that $G=HL$ if and only if $H$ acts transitively on $G/L$ and, in the same way, $L$ acts transitively on $G/H$. On the Lie algebra level, if $G=HL$ then $\mathfrak{g}=\mathfrak{h}+\mathfrak{l}$. The triple $(\mathfrak{g},\mathfrak{h},\mathfrak{l})$ is also called the Onishchik triple. In the case of a simple Lie algebra $\mathfrak{g}$ all possible Onishchik triples were classified in \cite{o1}, Table 2. In this table one finds also the intersections $\mathfrak{h}\cap\mathfrak{l}$. As a result, one can write down all standard triples $(\mathfrak{g},\mathfrak{h},\mathfrak{l})$. Note that in \cite{o1} the subalgebras $\mathfrak{h}$ and $\mathfrak{l}$ are denoted by $\mathfrak{g}'$ and $\mathfrak{g''}$. These observations yield the following.
\begin{proposition}\label{prop:onish-clas} The standard triples are contained in  the following table.

\begin{center}
 \begin{table}[ht]
 %\resizebox{\textwidth}{!}{
 \centering
 \begin{tabular}{| c | c | c | c | }
   \hline
	\textbf{$\mathfrak{g}$} & \textbf{$\mathfrak{h}$} & \textbf{$\mathfrak{l}$} & \textbf{$\mathfrak{h}\cap \mathfrak{l}$} \\
	 \hline
$\mathfrak{su}(2p,2)$ & $\mathfrak{sp}(p,1)$ & $\mathfrak{su}(2p,1)$ & $\mathfrak{sp}(p)$ \\
   \hline
$\mathfrak{su}(2,2(n-1)$ & $\mathfrak{sp}(1,n-1)$ & $\mathfrak{su}(1,2(n-1)$ & $\mathfrak{sp(n-1)}$ \\
   \hline
$\mathfrak{so}(3,4)$ & $\mathfrak{g}_{2(2)}$ & $\mathfrak{so}(1,4)$ & $ \mathfrak{su}(2)$ \\
   \hline
$\mathfrak{so}(2p,2)$ & $\mathfrak{so}(2p,1)$ & $\mathfrak{su}(p,1)$ & $\mathfrak{su}(p)$ \\
   \hline
$\mathfrak{so}(2,2(n-1))$ & $\mathfrak{so}(1,2(n-1))$ & $\mathfrak{su}(1,n-1)$ & $\mathfrak{su}(n-p)$ \\
   \hline
$\mathfrak{so}(4(n-1),4)$ & $\mathfrak{so}(4(n-1),3)$ & $\mathfrak{sp}(n-1,1)$ & $\mathfrak{sp}(p)$ \\
   \hline
$\mathfrak{so}(4(n-1),4)$ & $\mathfrak{so}(4(n-1),3)$ & $\mathfrak{sp}(n-1,1)$ & $ \mathfrak{sp}(n-1)\times\mathfrak{sp}(1)$ \\
   \hline
$\mathfrak{so}(8,8)$ & $\mathfrak{so}(7,8)$ & $\mathfrak{so}(1,8)$ & $\mathfrak{so}(7)$ \\
   \hline
$\mathfrak{so}(4,4)$ & $\mathfrak{so}(3,4)$ & $\mathfrak{so}(1,4)$ & $\mathfrak{so}(3)$ \\
   \hline
$\mathfrak{so}(4,4)$ & $\mathfrak{so}(3,4)$ & $\mathfrak{so}(1,4)\times \mathfrak{so}(3)$ & $\mathfrak{so}(3)\times\mathfrak{so}(3)$ \\
   \hline 
 \end{tabular}
 %}
 \caption{
 The standard triples
 }
 \label{t1}
 \end{table}
\end{center}
\end{proposition}
\newpage
\subsection{Einstein metrics on dual homogeneous spaces}
Let $G/H$ be a homogeneous space of semisimple non-compact Lie group $G$ and a Lie subgroup $H$ reductive in $G$. Denote by $G_{u}/H_{u}$ the compact dual of $G/H$. Recall the following theorem of I. Kath \cite{kat}.
\begin{theorem}[\cite{kat}, Corollary 4.1 and 4.2]\label{thm:kath}
The $G$-invariant pseudo-  Riemannian Einstein metrics on $G/H$ are in bijective correspondence with $G_{u}$-invariant Riemannian Einstein metrics on $G_{u}/H_{u}.$ 
\label{kkaa}
\end{theorem}
Let us make the following observation. Assume that a reductive subgroup $L$ acts properly and co-compactly on $G/H$. Then $H$ acts properly and co-compactly on $G/L$. This is easy to see or one can consult \cite{kob}. Therefore, any standard triple $(G,H,L)$ yields {\it two types} of compact Clifford-Klein forms: $\Gamma_1\backslash G/H$ and $L\backslash G/\Gamma_2$ with co-compact lattices $\Gamma_1\subset L$ and $\Gamma_2\subset H$. For the convenience of the reference we will call the homogeneous spaces obtained from the standard triples as {\it Kobayashi-Onishchik spaces}, or {\it KO-spaces} for short.

\begin{proposition}\label{prop:ko-spaces}The following table contains all spaces dual to KO-spaces:
\begin{center}
 \begin{table}[ht]
 \resizebox{\textwidth}{!}{%
 \centering
 {\footnotesize
 \begin{tabular}{| c | c | c | c | c |}
   \hline
   \multicolumn{5}{|c|}{\textbf{Table 2.} \textbf{\textit{Compact duals of KO-spaces ($n\geq 1$)}}} \\
   \hline                        
    & KO-spaces $G/H$ & $G/H'$ ($H'\subset H$-without compact factors) & Compact dual $G_{u}/H_{u}$ & Compact dual $G_{u}/H'_{u}$\\
   \hline
   1 & $SU(2,2n)/Sp(1,n)$  & $SU(2,2n)/Sp(1,n)$  & $SU(2n+2)/Sp(n+1)$ &  $SU(2n+2)/Sp(n+1)$ \\
   \hline
   2 & $SU(2,2n)/U(1,2n)$  & $SU(2,2n)/SU(1,2n)$  & $SU(2n+2)/U(2n+1)$ & $SU(2n+2)/SU(2n+1)$ \\
   \hline
   3 & $SO(2,2n)/U(1,n)$  & $SO(2,2n)/SU(1,n)$  & $SO(2n+2)/U(n+1)$ & $SO(2n+2)/SU(n+1)$ \\
   \hline
   4 & $SO(2,2n)/SO(1,2n)$  & $SO(2,2n)/SO(1,2n)$  & $SO(2n+2)/SO(2n+1)$ & $SO(2n+2)/SO(2n+1)$ \\
   \footnotesize $n\geq 2$  & & & & \\ 
   \hline
	 5 & $SO(4,4n)/SO(3,4n)$  & $SO(4,4n)/SO(3,4n)$  & $SO(4n+4)/SO(4n+3)$ & $SO(4n+4)/SO(4n+3)$ \\
   \hline
   6 & $SO(4,4n)/Sp(1,n)Sp(1)$  & $SO(4,4n)/Sp(1,n)$  & $SO(4n+4)/Sp(n+1)Sp(1)$ & $SO(4n+4)/Sp(n+1)$ \\
   \hline
	 7 &  $SO(4,4)/SO(4,1) \times SO(3)$ & $SO(4,4)/SO(4,1)$ & $SO(8)/SO(5) \times SO(3)$ & $SO(8)/SO(5)$ \\
	 \hline
   8 &  $SO(4,4)/Spin(4,3)$  & $SO(4,4)/Spin(4,3)$  & $SO(8)/Spin(7)$ & $SO(8)/Spin(7)$ \\
   \hline  
	 9 &  $SO(4,3)/SO(4,1) \times SO(2)$  & $SO(4,3)/SO(4,1)$ & $SO(7)/SO(5)\times SO(2)$ & $SO(7)/SO(5)$ \\
   \hline 
   10 & $SO(4,3)/G_{2(2)}$  & $SO(4,3)/G_{2(2)}$ & $SO(7)/G_{2}$ & $SO(7)/G_{2}$ \\
   \hline 
	 11 &  $SO(8,8)/SO(7,8)$  & $SO(8,8)/SO(7,8)$  & $SO(16)/SO(15)$ & $SO(16)/SO(15)$ \\
   \hline 
	 12 &  $SO(8,8)/Spin(1,8)$  & $SO(8,8)/Spin(1,8)$  & $SO(16)/Spin(9)$ & $SO(16)/Spin(9)$ \\
   \hline 
 \end{tabular}}
 }
 \caption{
 KO-spaces
 }
 \label{t2}
 \end{table}
\end{center}
\end{proposition}
\newpage 
\subsection{Inspection of Table 2}
Recall that Theorem \ref{kkaa} shows that we need to analyze compact homogeneous spaces from Table 2 (which we call Kobayashi-Onishchik triples).  
We begin with the following results from \cite{wz} (Corollary 2.3), \cite{be} (Theorem 7.44) and \cite{W}.
\begin{proposition}\label{prop:wz} Let $H$ be any closed subgroup of a compact connected Lie group $G$. If $\mathfrak{h}$ is maximal in $\mathfrak{g}$, then $G/H$ admits a $G$-invariant Einstein metric.
\end{proposition}
\begin{proposition}\label{prop:isotropy} If $G/H$ is isotropy irreducible, then it admits a unique (up to a scalar)  $G$-invariant metric, which is Einstein.
\end{proposition}
Using Propositions \ref{prop:wz} and \ref{prop:isotropy} and Table 7.107 in \cite{be}, we remove from Table 2 all isotropy irreducible homogeneous spaces, all homogeneous spaces $G/H$ with $\mathfrak{h}$ maximal in $\mathfrak{g}$ and all standard spheres (because all of them admit Einstein metrics, and we don't need to consider them).  
 Next, we write down the remaining triples which are
\begin{itemize}
\item spheres $SU(2n+2)/SU(2n+1)$, \ \ \ \ \ $SU(2n+1)\subset U(2n+1),$
\item homogeneous spaces $SO(2n+2)/SU(n+1),$ determined by a standard inclusion  $SU(n+1) \subset U(n+1),$
\item Stiefel manifolds $SO(8)/SO(5),$ and $SO(7)/SO(5)$
\item  homogeneous manifolds $SO(4n)/Sp(n)$ and $SO(4n)/Sp(n) Sp(1)$.
\end{itemize}
\section{Invariant Einstein metrics with additional symmetries}\label{subsec:symmetry}
In our proof we will use a method \cite{aa} of constructing invariant Einstein metrics on compact homogeneous spaces $G/H$ using additional symmetries. Here we briefly explain it. Let $G/H$ be a compact homogeneous space. Consider the decomposition $\mathfrak{g}=\mathfrak{h}+\mathfrak{m},\mathfrak{m}\cong T_HG/H$.  The normalizer $N_G(H)$ acts on $G/H$ by $(a,gH)\rightarrow ga^{-1}H$. If $a\in N_G(H)$ is fixed, one obtains a diffeomorphism $\varphi_a: G/H\rightarrow G/H,\,\varphi_a(gH)=ga^{-1}H$. One easily checks that if $\langle\cdot,\cdot\rangle$ is an $Ad(H)$-invariant inner product determining a $G$-invariant Riemannian metric $\rho$ on $G/H$, then $\varphi_a$ is an isometry if and only if $\operatorname{Ad}(a)|_{\mathfrak{m}}$ preserves $\langle\cdot,\cdot\rangle$. Assume that there exists a closed subgroup $S\subset G$ which centralizes $H$.
Let $Z=H\cap S$. Since $S$ centralizes $H$, $Z$ is central in $H$ and $S$. 
 Then   $\tilde G=G\times (S/Z)$ acts on $G/H$ by $(a,[s])\cdot gH=ags^{-1}H$. The isotropy group $\tilde H$ of this action is isomorphic to $K=HS\subset G$, the isomorphism is given by $([s],sh)\rightarrow sh$.
Let $\mathcal{M}^{G,K}$ denote the set of $G$-invariant metrics on $G/H$ which correspond to the inner products on $\mathfrak{m}$ which are $\operatorname{Ad}(K)$ Then the subgroup -invariant. Therefore 
$$\mathcal{M}^{G,K}\subset\mathcal{M}^G.$$
Thus, the set $\mathcal{M}^{G,K}$ consists of $G$-invariant metrics with additional symmetry (they are determined by $\operatorname{Ad}(K)$-invariant inner products, where $K=H S$). Also, one may note that $\varphi(a)$ are isometries for such metrics, the $\tilde G=G\times S/Z$ acts on $G/H$ by isometries, and one may identify $\mathcal{M}^{\tilde G}$ with $\mathcal{M}^{G,K}$. 
\begin{definition} Metrics in $\mathcal{M}^{G,K}$ will be called $\operatorname{Ad}(K)$-invariant metrics on $G/H$. If $K\cong H S$ with $S\subset Z_G(H)$ and with finite $H\cap S$, we will say that $K$ is of a centralizer type.
\end{definition}
Now, consider only $G$-invariant Riemannian metrics $\rho$ on $G/H$ which are also $\operatorname{Ad}(K)$-invariant (thus, $\rho\in\mathcal{M}^{G,K}$). 

\begin{proposition}[\cite{aa}, Proposition 3]
Let $G$ be compact, connected  and semisimple. Assume that we are given inclusions $H\subset K\subset G$ and that $K$ is of a centralizer type. The $\operatorname{Ad}(K)$-invariant metric of volume one on $G/H$ is Einstein if and only if it is a critical point of the scalar curvature functional $S$ restricted to $\mathcal{M}^{G,K}_1$.
\label{propk}
\end{proposition}
\textsl{Proof.} The proof follows since one can identify $\mathcal{M}^{\tilde G}$ with $\mathcal{M}^{G,K}$ (as in subsection \ref{subsec:symmetry}) and then use the general fact that invariant Einstein metrics on homogeneous spaces of compact Lie groups are critical points of the restriction of the scalar curvature functional \cite{be}, Theorem 4.23.

\hfill$\square$

\noindent \textcolor{blue}{It shoud be noted that the method of constructing Einstein metrics with additional symmetries is explained in more detail in \cite{st}.} 
\section{Proof of Theorem \ref{thm:main}} 
We complete the proof as a case-by-case analysis. It appears that in all cases there are $G$-invariant Einstein metrics, although there is no complete classification. We begin with the known cases. 
\begin{itemize}
\item {\bf The case of spheres $S^{4n+3}=SU(2n+2)/SU(2n+1)$}. This is known to admit (non-standard) invariant Einstein metric \cite{J}, \cite{z}.
\item  {\bf The case $SO(2n+2)/SU(n+1)$.} We consider this case as a homogeneous fiber bundle
$$S^1\rightarrow SO(2n)/SU(n)\rightarrow SO(2n)/U(n)$$
and use the result of Jensen \cite{J} to state that the total space admits an invariant Einstein metric.
\item {\bf The case $SO(8)/SO(5)$}. This is the Stiefel manifold of the form $$SO(l+k)/SO(k)$$
 with $k=3$. Referring to \cite{J} again we deduce that there is at least 2 $G$-invariant Einstein metrics.
\item {\bf The space $SO(8)/SO(5)\times SO(2)$}. This is a {\it generalized Wallach space} of the form $SO(k+l+m)/SO(k)\times SO(l)\times SO(m)$ with $m=1$. We refer to \cite{cn} \textcolor{blue}{for the existence of Einstein metrics.}
\item {\bf The case $SO(7)/SO(5)$}. This is the  Stiefel manifold which admits a $G$-invariant Einstein metric by \cite{bh}.
\end{itemize}

Now we consider  {\bf the cases  $SO(4n)/Sp(n)Sp(1)$ and $SO(4n)/Sp(n)$.}
\vskip6pt
To our best knowledge they are not covered by the previous work. First notice that $SO(4n)/Sp(n)Sp(1)$ is isotropy irreducible (for example see the list on page 185 in \cite{wz}) and thus it admits a unique (up to a constant) invariant Einstein metric. We also have the following.
\begin{theorem}
The space $$SO(4n)/Sp(n)$$ admits at least two different invariant Einstein metrics.
\end{theorem}
To prove this theorem we will use the method described in Section \ref{subsec:symmetry}. 
 
%We consider only $G$-invariant Riemannian metrics $\rho$ on $G/H$ which are also $\operatorname{Ad}(K)$-invariant (thus, $\rho\in\mathcal{M}^{G,K}$). 

We  make the following straightforward observation.
\begin{proposition} \label{prop:adk-isotropy}Let $G$ be compact, connected  and semisimple. Assume that we are given the inclusions $H\subset K\subset G$ and the corresponding decompositions of Lie algebras
$$\mathfrak{g}=\mathfrak{h}+\mathfrak{m_1}+\mathfrak{m}_2,\mathfrak{k}=\mathfrak{h}+\mathfrak{m}_1.$$
If $\mathfrak{m}_1+\mathfrak{m}_2$ is a decomposition into isotropy irreducible $\operatorname{Ad}(K)$-modules, then the inner products
$$\langle\cdot,\cdot\rangle=x_1B|_{\mathfrak{m}_1}+x_2B|_{\mathfrak{m}_2},\,x_i>0,i=1,2 \eqno (6)$$
determine all $\operatorname{Ad}K$-invariant Riemannian metrics on $G/H$.
\end{proposition}   
We get the following.
\begin{corollary}\label{cor:method} \textcolor{blue}{Let $G$ be compact connected and semisimple. Assume that we are given inclusions $H\subset K\subset G$ and that $K$ is of centralizer type such that
\begin{itemize}
	\item  $G/K$ is isotropy irreducible,
	\item the action of K on $T(K/H)$ is irreducible.
\end{itemize}
Then the $\operatorname{Ad}K$-invariant Einstein metrics on $G/H$ are of the form
$$\langle\cdot,\cdot\rangle=x_1B|_{\mathfrak{m}_1}+x_2B|_{\mathfrak{m}_2},$$
where $x_{1}, x_{2}$ are given by formula $(5).$}
\end{corollary}
\textsl{Proof.} Under the assumptions of the Corollary, $\mathfrak{m}_i, i=1,2$ are irreducible $\operatorname{Ad}(K)$-modules. Therefore, by Proposition \ref{prop:adk-isotropy} all $\operatorname{Ad}(K)$-invariant metrics on $G/H$ are determined by formula $(6)$. We know by Proposition \ref{propk} that the critical points of $S|_{\mathcal{M}^{G,K}}$ are Einstein. On the other hand, we know by Lemma \ref{lagrange} that these critical points are calculated by formula $(5)$. \hfill$\square$

Thus we are looking for the invariant Einstein metrics on $G/H$ as follows. Set $G=SO(4n),$ $n\geq 2,$ $K=Sp(n)Sp(1)$ and $H=Sp(n), L=Sp(1)$. We check, if there are Einstein metrics in $\mathcal{M}^{G,K}_1$. By Corollary \ref{cor:method}  it is sufficient to show that the polynomial (5) has  real roots. 
\begin{lemma}\label{lemma:Killing-id} Let $\mathfrak{g}$ be a compact semisimple Lie algebra of dimension $\dim\mathfrak{g}=c$ and $B$ the negative of the Killing form of $\mathfrak{g}$.  For any   $B$-orthonormal basis $\{ h_{i} \}$ the following equality holds
$$\sum_{j,k=1}^{c} (B([h_{i},h_{j}],h_{k}))^{2}=\sum_{j=1}^{c}B([h_{i},h_{j}],[h_{i},h_{j}])=B(h_{i},h_{i})=1.$$
\end{lemma}
\textsl{Proof of Lemma \ref{lemma:Killing-id}}.
This is a direct calculation using the  properties of the Killing form and the expression
$$B(h_i,h_j)=\sum_{k,l=1}^cc_{il}^kc_{jk}^l$$
(compare \cite{aa}, the proof of Lemma 4). Here $c_{ij}^k$ denote the structural constants of $\mathfrak{g}$.

\hfill$\square$

\begin{lemma}\label{lemma:poly}
If $G=SO(4n),$ $n\geq 2,$ $K=Sp(n)Sp(1)$ and $H=Sp(n), L=Sp(1)$ then the coefficients of polynomial $(5)$ are given by the formulas
$$
[1,1,1]=\frac{3n+3}{4n-2},\, [1,2,2]=\frac{9n-9}{4n-2},$$
 $$[2,2,2]=3(2n^{2}-n-1)\frac{n^{2}-n-2}{2n^{2}-n}.\eqno (7)
 $$
\end{lemma}
\textsl{Proof of  Lemma \ref{lemma:poly}}.
Recall that  $B$ denotes the negative of the Killing form of $\mathfrak{g}$, which is positively definite. For a subalgebra $\mathfrak{k}\subset\mathfrak{g}$ we denote by $B_{\mathfrak{k}}$ the negative of the Killing form of $\mathfrak{k}$, while $B|_{\mathfrak{k}}$ is the restriction of $B$ onto $\mathfrak{k}$.  
Denote by $c^{\ast}_{2}$ the Casimir constant of $G/K$ with respect to $B |_{\mathfrak{k}}.$ Note that in under the conditions of the lemma $d_{1}=3,$ $d_{2}=6n^2-3n-3.$ Also $\alpha B|_{\mathfrak{k}}=B_{\mathfrak{k}},$ $\alpha = \frac{n+1}{4n-2}$ and $c^{\ast}_{2}=\frac{n^2+2}{4n^2-2n}$ (\cite{wz1}, Table V).

Because $G/K$ is isotropy irreducible and $\mathfrak{g}=\mathfrak{k}+\mathfrak{m}_{2}$ we may compute $[2,2,2]$ using Lemma \ref{lemma:casimir}. Thus, 
$$[2,2,2]=d_{2}(1-2c^{\ast}_{2})=3(2n^{2}-n-1)\frac{n^2-n-2}{2n^2-n}.$$

To compute [1,1,1] we use Lemma \ref{lemma:Killing-id}. Take the $B$-orthonormal  basis of $\mathfrak{m}_{1}=\mathfrak{sp}(1),$ $\{ e_{\mu} \} \cap \mathfrak{m}_{1} := \{ e_{1},e_{2},e_{3}\} .$ Since $\mathfrak{m}_{1}$ is an ideal in $\mathfrak{k}$  (see \cite{wd2}, Theorem 11)  $B_{s}:=\alpha B|_{\mathfrak{m}_{1}}$ is the negative of the Killing form of $\mathfrak{m}_{1}$ and $\{ v_{i}:=\frac{e_{i}}{\sqrt{\alpha}} \} $ is an orthonormal basis of $\mathfrak{m}_{1}.$ We have
$$[1,1,1]=\sum_{i=1}^{3} (\sum_{j,k=1}^{3} (B([e_{i},e_{j}],e_{k}))^{2})=\sum_{i=1}^{3} (\sum_{j=1}^{3}B([e_{i},e_{j}],[e_{i},e_{j}]))=$$
$$=\frac{1}{\alpha}\sum_{i=1}^{3}(\sum_{j=1}^{3}B_{s}([e_{i},e_{j}],[e_{i},e_{j}]))=\alpha \sum_{i=1}^{3} (\sum_{j=1}^{3} B_{s}([v_{i},v_{j}],[v_{i},v_{j}]))=$$
$$=\alpha \sum_{i=1}^{3} B_{s}(v_{i},v_{i})=3\alpha ,$$
because of Lemma \ref{lemma:Killing-id} and $\dim\mathfrak{m}_{1}=3.$

To calculate $[1,2,2]$, put $\mathfrak{m}_{0}:=\mathfrak{h}$ and extend the basis $\{ e_{\mu} \}$ with an orthonormal (with respect to  $B$) basis of $\mathfrak{m}_{0}$ so that if $e_{\mu}\in \mathfrak{m}_{i}$ and $e_{\beta}\in \mathfrak{m}_{j}$ and $i<j$ then $\mu < \beta .$ We have the following equalities
$$3= \sum_{e_{\mu} \in \mathfrak{m}_{1}} B(e_{\mu}, e_{\mu})= \sum_{e_{\mu} \in \mathfrak{m}_{1}} \sum_{e_{\beta},e_{\gamma} \in \mathfrak{m}_{0}+\mathfrak{m}_{1}+\mathfrak{m}_{2}} (B([e_{\mu},e_{\beta}],e_{\gamma}))^{2}.$$
Since $[\mathfrak{m}_{0},\mathfrak{m}_{1}]=0$ one can continue as follows
$$\sum_{e_{\mu} \in \mathfrak{m}_{1}} \sum_{e_{\beta},e_{\gamma} \in \mathfrak{m}_{0}+\mathfrak{m}_{1}+\mathfrak{m}_{2}} (B([e_{\mu},e_{\beta}],e_{\gamma}))^{2}
=\sum_{e_{\mu} \in \mathfrak{m}_{1}} \sum_{e_{\beta},e_{\gamma} \in \mathfrak{m}_{1}+\mathfrak{m}_{2}} (B([e_{\mu},e_{\beta}],e_{\gamma}))^{2}$$
$$=[1,1,1]+[1,1,2]+[1,2,1]+[1,2,2].$$
 We know that $\mathfrak{k}=\mathfrak{m}_{0}+\mathfrak{m}_{1}$ is a subalgebra, therefore $[1,1,2]=[1,2,1]=0$. Finally
$$[1,2,2]=3-[1,1,1]$$
and the proof of the Lemma follows.

\hfill$\square$

Now we complete the proof of Theorem \ref{thm:main} by a straightforward calculation of   polynomial (5)  which has the final form
$$3(n-1)(2n^{2}-n+1)t^{2}-2(2n^{2}-n-1)\frac{3n^{2}-n+2}{n}t+(n+1)(2n^{2}-n-1).$$
One can compute that it has two real roots. The proof is complete.

 \hfill$\square$

Department of Mathematics and Computer Science, 

University of Warmia and Mazury,

S\l\/oneczna 54, 10-710, Olsztyn, Poland

mabo@matman.uwm.edu.pl (MB),

tralle@matman.uwm.edu.pl  (AT).

\end{document}